\documentclass{amsart} % in was article
\usepackage{amssymb}
\usepackage{amsmath}
\usepackage{amsfonts}

\begin{document}
%%%%%%%%%%%%%%%%%%%%%%%%%%%%%%%%%%%%%%%%%%%%%%%%%%%%%%%%%%%%%%%%%%%%%%
%	spaces for your own definitions follows
%%%%%%%%%%%%%%%%%%%%%%%%%%%%%%%%%%%%%%%%%%%%%%%%%%%%%%%%%%%%%%%%%%%%%%
\newtheorem{theo}{Theorem}[section]
\newtheorem{prop}[theo]{Proposition}
\newtheorem{lemma}[theo]{Lemma}
\newtheorem{exam}[theo]{Example}
\newtheorem{coro}[theo]{Corollary}
\theoremstyle{definition}
\newtheorem{defi}[theo]{Definition}
\newtheorem{rem}[theo]{Remark}

%\renewcommand{\theequation}{\mbox{\arabic{section}.\arabic{equation}}}

%letters - added these
\newcommand{\Bb}{{\bf B}}
\newcommand{\Nb}{{\bf N}}
\newcommand{\Qb}{{\bf Q}}
\newcommand{\Rb}{{\bf R}}
\newcommand{\Zb}{{\bf Z}}
\newcommand{\Ac}{{\mathcal A}}
\newcommand{\Bc}{{\mathcal B}}
\newcommand{\Cc}{{\mathcal C}}
\newcommand{\Dc}{{\mathcal D}}
\newcommand{\Fc}{{\mathcal F}}
\newcommand{\Ic}{{\mathcal I}}
\newcommand{\Jc}{{\mathcal J}}
\newcommand{\Lc}{{\mathcal L}}
\newcommand{\Oc}{{\mathcal O}}
\newcommand{\Pc}{{\mathcal P}}
\newcommand{\Sc}{{\mathcal S}}
\newcommand{\Tc}{{\mathcal T}}
\newcommand{\Uc}{{\mathcal U}}
\newcommand{\Vc}{{\mathcal V}}

\newcommand{\ax}{{\rm ax}}
\newcommand{\Acc}{{\rm Acc}}
\newcommand{\Act}{{\rm Act}}
\newcommand{\ded}{{\rm ded}}
\newcommand{\Gm}{{$\Gamma_0$}}
\newcommand{\ID}{{${\rm ID}_1^i(\Oc)$}}
\newcommand{\PAP}{{${\rm PA}(P)$}}
\newcommand{\ACA}{{${\rm ACA}^i$}}
\newcommand{\RefP}{{${\rm Ref}^*({\rm PA}(P))$}}
\newcommand{\RefS}{{${\rm Ref}^*({\rm S}(P))$}}
\newcommand{\Rfn}{{\rm Rfn}}
\newcommand{\tar}{{\rm Tarski}}
\newcommand{\UNFA}{{${\mathcal U}({\rm NFA})$}}

\author{Nik Weaver}

\title [Mathematical conceptualism]
       {Mathematical conceptualism}

\address {Department of Mathematics\\
          Washington University in Saint Louis\\
          Saint Louis, MO 63130}

\email {nweaver@math.wustl.edu}

\date{\em September 12, 2005}

\maketitle

%%%%%%%%%%%%%%%%%%%%%%%%%%%%%%%%%%%%%%%%%%%%%%%%%%%%%%%%%%%%%%%%%%%%%%%
%	Please insert the article body now
%%%%%%%%%%%%%%%%%%%%%%%%%%%%%%%%%%%%%%%%%%%%%%%%%%%%%%%%%%%%%%%%%%%%%%%

\begin{quote}
We must remember that in Hilbert's time all mathematicians were excited
about the foundations of mathematics. Intense controversy centered around
the problem of the legitimacy of abstract objects $\ldots$ The contrast
with today's foggy atmosphere of intellectual exhaustion and
compartmentalization could not be more striking.

--- Stephen Simpson (\cite{Sim1}, p.\ 350)
\end{quote}

From the standpoint of mainstream mathematics, the great foundational
debates of the early twentieth century were decisively settled in favor
of Cantorian set theory, as formalized in the system ZFC (Zermelo-Fraenkel
set theory including the axiom of choice). Although basic foundational
questions have never entirely disappeared, it seems fair to say that they
have retreated to the periphery of mathematical practice. Sporadic
alternative proposals like topos theory or Errett Bishop's constructivism
have never attracted a substantial mainstream following, and Cantor's
universe is generally acknowledged as the arena in which modern mathematics
takes place.

Despite this history, I believe a strong case can be made for abandoning
Cantorian set theory as the de facto foundational standard, in favor of an
alternative stance which I call {\it mathematical conceptualism}. An analogy
can be made with the nineteenth-century transition from the naive use of
infinitesimals to the epsilon-delta method. In both cases the defining feature
of the transition is the rejection of a class of ill-defined metaphysical
entities in favor of a more concrete and detailed picture. In both
cases, too, although from a philosophical point of view the transition is
rather drastic, mainstream mathematical practice is not fundamentally
altered (though it is made more rigorous). At a deeper level, there is in
both cases a sense in which the old view can be recovered under the new one.
The analogy even extends a bit further; I will return to it later.

Mathematical conceptualism is a modern version of an older philosophy called
predicativism. Although predicativism is practically unknown to working
mathematicians today, it has an impressive pedigree: it was originally
formulated by
Henri Poincar\'e and Bertrand Russell, and Hermann Weyl made a major early
contribution. Moreover, its central principle, that sets must be ``built
up from below'', seems so inoffensive that many mathematicians might be
surprised to learn that they routinely violate it. So predicativism's rapid
decline into obscurity may seem puzzling. The principal reason for this
decline was undoubtedly that it was initially seen as being too weak to support
mainstream mathematics. Although we now know that this is not true, it seems
for instance to have inhibited Russell even from fully accepting his own
ideas. Poincar\'e and Weyl, on the other hand, held less fixed philosophical
views and so predicativism was left without a true believer to defend it
such as intuitionism had in L.\ E.\ J.\ Brouwer.

But even if it had had such a single-minded advocate, and even if its true
strength had been recognized earlier, predicativism would probably still have
been unsustainable because it conflicted with the basic twentieth-century
trend towards abstraction. However, I believe that trend has basically lost
its force and abstraction for its own sake is no longer an attractive idea
to many working mathematicians. On the contrary, we have now reached a point
where we have a fairly clear idea of just what portion of the Cantorian
universe is relevant to mainstream mathematics, and the fact that this region
is, with remarkable accuracy, precisely the portion a conceptualist would
recognize as legitimate may be cited as a point in conceptualism's favor
(see Section \ref{sect5}).

A more perplexing factor contributing to the near total eclipse of
predicativism in recent years has been its treatment by workers in
mathematical foundations, who might have been expected to retain real
interest in it but instead for the past forty years have uncritically
endorsed an incorrect analysis which made it appear trivially false. This
may seem like a strong charge but I am afraid it is accurate \cite{Wea2}.

This paper is non-technical and is addressed to a general mathematical
and philosophical audience. In two companion papers I discuss in some
detail the practical matter of how one is to actually do mathematics
conceptualistically \cite{Wea1} and the proof-theoretic strength of
conceptualistic theories \cite{Wea2}. The first of these should be
accessible to mathematicians and the second to logicians.

\section{Anomalies in set theory}\label{sect1}

\begin{quote}
Set theory, as formalized in ZF, provides an extremely powerful and
elegant way to organize existing mathematics. It is not without its
drawbacks, nevertheless. While it is too weak to decide some questions
(like the continuum hypothesis) which seem meaningful (even important),
it is in some ways too strong.

--- Jon Barwise (\cite{Bar}, pp.\ 7-8)
\end{quote}

Before I discuss conceptualism, I want to give the reader some idea of
what I think is wrong with Cantorian set theory. The reader, of course,
may already have ideas of his own on this question --- I am unsure of
the extent to which most working mathematicians are really ideologically
committed to the classical picture.

I call the points I raise below ``anomalies'' because I do not claim
any of them rises to the level of discrediting classical set theory
outright. However, they at least ought to make one pause. It will become
apparent that the conceptualist stance neatly answers all of them.
\medskip

\noindent {\bf Anomaly \# 1}: {\it Undecidability of the continuum
hypothesis.}

According to G\"odel's incompleteness theorem (as improved by Rosser), no
consistent formal system strong enough to support elementary number theory is
complete, i.e., there necessarily exist statements whose truth value cannot be
determined within the system. However, the actual examples of such statements
offered by G\"odel are number-theoretically rather complicated and quite
remote from mainstream mathematical practice.
Harvey Friedman has made a valiant and technically brilliant effort to produce
simple and compelling number-theoretic statements that are undecidable in
standard set theory (e.g., see \cite{Fri}), but the verdict on this work
still seems to be out. It is certainly very interesting in that he produces
elementary number-theoretic assertions that are formally not terribly
complicated and have no obviously transcendent character but are in fact
equivalent to the consistency with ZFC of the existence of various large
cardinals. However, their relevance to mainstream mathematics is not so
clear. (Friedman has also claimed that some manifestly mainstream results
are independent of predicative reasoning, but this argument depends on the
faulty analysis of predicativism that I mentioned earlier.)

What is more disturbing about Cantorian set theory is that its standard
formalizations are known to be incapable of answering questions of seemingly
basic importance, most prominently the continuum hypothesis. Now this in
itself is obviously not enough to actually falsify the Cantorian view, but it
is certainly unsettling. As Paul Cohen has written, ``This state of affairs
regarding a classical and presumably well-posed problem must certainly appear
rather unsatisfactory to the average mathematician'' (\cite{Coh}, p.\ 1). I
would add that the problem is not only classical and presumably well-posed,
but also apparently quite fundamental.

Recently Hugh Woodin has proposed that the continuum hypothesis be
settled negatively using a new axiom related to projective determinacy
\cite{Woo1, Woo2}. Like Friedman's work mentioned above, Woodin's work
is clearly of great technical interest, but the claims he makes for it
seem too strong since the new axiom is rather complicated and the
arguments for accepting it are indirect to say the least. (Woodin's actual
views on the nature of mathematical truth are somewhat unusual; he has
elsewhere suggested, or seemed to suggest, that elementary number theory
has no canonical model and that if the twin primes conjecture is undecidable
in ZFC then it may have no definite truth value (\cite{Woo0}, p.\ 34). This
apparently indicates a disbelief in the natural numbers coupled with a
belief in the axioms of set theory, an odd combination.)
\medskip

\noindent {\bf Anomaly \# 2}: {\it Poor fit with mathematical practice.}

Classical set theory presents us with a picture of an incredibly vast universe.
This universe not only contains the set of natural numbers $\Nb$ and its power
set $\Pc(\Nb)$ --- the set of all sets of natural numbers --- it also contains
$\Pc(\Pc(\Nb))$, $\Pc(\Pc(\Pc(\Nb)))$,
etc. This process can be iterated uncountably many times. It can be iterated
$\aleph_\omega = \sup_n \aleph_n$ times or $\aleph_{\aleph_\omega}$ times.
It could be iterated $2^{\aleph_0}$ times, if we knew the value of $\alpha$
for which $2^{\aleph_0} = \aleph_\alpha$. {\it Yet virtually all important
objects in mainstream mathematics are either countable or separable.} This
is not because the uncountable/nonseparable case has not yet been sufficiently
studied but rather because on mainstream questions it tends to be either
pathological or undecidable. (I have had personal experience with both
cases \cite{AW, Wea0}.)

Apparently, out of the entire unimaginably vast set-theoretic universe,
only the bottommost layer is interesting from the point of view of
mainstream mathematics. Something seems wrong.

Again, this is not an utterly damning criticism of Cantorian set theory, but
it does indicate a basic disconnection between set theory and mainstream
mathematics. From the former point of view the latter is quite trivial and
from the latter point of view the former is almost completely irrelevant.
This characterization could be disputed, but I definitely think that Barwise's
comment quoted above that set theory ``is in some ways too strong'' is an
understatement from the point of view of the core mathematician.
\medskip

\noindent {\bf Anomaly \# 3}: {\it Skolem's paradox.}

According to the L\"owenheim-Skolem theorem, any infinite model of a countable
family of axioms has a countable submodel. This can be applied to the
standard axioms of set theory to infer the existence of a countable universe
in which the ZFC axioms hold. In other words, if we restrict our attention
to just those sets which belong to a certain countable set $M$, and ignore all
other sets, we find that the standard axioms of set theory hold. (This is
slightly subtle; perhaps the simplest accurate statement is that if
inaccessible cardinals exist then there is a countable transitive set
which models ZFC.)

Skolem's paradox is that since $M$ satisfies the ZFC axioms and we can
prove in ZFC that there exist uncountable sets, $M$ must therefore contain
uncountable sets, yet $M$ itself is countable. This puzzle is easily
resolved (and indeed the resolution was already offered by Skolem): $M$
``thinks'' that some of its elements are uncountable because it does not
contain a bijection between any of these sets and the natural numbers.
However, this does not prevent such bijections from existing outside $M$.

The real import of Skolem's paradox is the fact that $M$ mimics the actual
universe of sets so exactly that the latter becomes, for any mathematical
purpose, completely superfluous. Moreover, it seems to place the actual
universe so far out of our reach as to call our knowledge of its existence
into question. ``All commentators agree that the existence of such models
shows that the `intended' interpretation, or, as some prefer to speak, the
`intuitive notion of a set', is not `captured' by the formal system. But if
{\it axioms} cannot capture the `intuitive notion of a set', what possibly
could?'' (\cite{Put}, p.\ 465)

In addition, the technique of forcing gives us a great deal of freedom to
convert a given countable model $M$ into other models $M'$, $M''$, $\ldots$
with various properties. So apparently one is actually faced with a wide
variety of countable models of the ZFC axioms, while only one (or a handful)
of them is supposed to reflect the true properties of the real universe. Yet
we evidently have no way to determine which is which.
\medskip

\noindent {\bf Anomaly \# 4}: {\it The classical paradoxes.}

The well-known set-theoretic paradoxes demonstrate that naive set theory is
inconsistent. Ideally we would like to know in exactly what way the intuitions
underlying naive set theory are erroneous, and we would then want to adopt
set-theoretic axioms which express our corrected intuition. Instead, one has
a sense that the ZFC axioms were chosen in a somewhat ad hoc manner in which
the goal was not to achieve philosophical coherence, but rather simply to
retain as much of naive set theory as one could without admitting the known
paradoxes. (``There is, to be sure, a certain justification for the axioms
in the fact that they go into evident propositions of naive set theory if in
them we take the word `set', which has no meaning in the axiomatization, in
the sense of Cantor. But what is omitted from naive set theory --- and to
circumvent the antinomies some omission is essential --- is absolutely
arbitrary'' (\cite{vN}, p.\ 396).)

The result is an incoherent system that is based on contradictory intuitions.
It is supposed to spring from a conception of an iteratively generated
universe which is built up in stages, and indeed {\it most} of the axioms
fit this view. Yet the power set axiom is incompatible with it because we do
not visualize building up power sets in a step-by-step fashion. Instead, ZFC
simply postulates that once a set is available its power set immediately
becomes available too. This incongruity is so severe that I am surprised it
has not attracted more comment. If we really believe that the set-theoretic
universe has to be built up piecemeal we surely cannot accept an axiom
according to which enormous new sets (enormous because there is a jump in
cardinality) simply nonconstructively appear. Where do they come from, if
they were not in some sense already there? Conversely, if we believe that
whenever a set exists then its power set does not need to be constructed,
but is indeed simply ``there'', then we must believe this of the power set
of the power set, the power set of the power set of the power set, and so on.
In fact, would it not be the case that the entire universe is simply ``there''?

At bottom this is a question about what brings the set-theoretic universe into
existence. Explanations of the iterative conception are completely vague on
this point. How does one reconcile the belief that sets may exist without
being explicitly constructible with the idea of a universe that is generated
in stages? If {\it we} are not responsible for constructing the universe,
then what causes it to appear? Is this a process which takes place over time?
Is it possible to imagine a dynamic process which does {\it not} take place
over time?

This leads into the question of what sets actually are, which I will address
in Section \ref{sect2}, so I will stop here with the one additional comment
that the most obvious way to justify the existence of a set that one is unable
to explicitly construct is in the context of a preexisting universe out of
which one could imagine globally selecting the elements of the desired set.
But the assumption of a well-defined preexisting universe of sets not only
flatly contradicts the iterative conception, it also appears to support the
construction of paradoxical sets (e.g., Russell's set) just as well as it
supports the construction of power sets.
\medskip

All four anomalies fit together: each, in one way or another, expresses
the idea that the Cantorian universe is too big. Its construction
incorporates a power set operation whose intuitive justification is
unclear and which makes the universe vastly larger than it would be
otherwise (anomaly \# 4). The great size of the resulting universe is
both philosophically dubious (anomaly \# 3) and mathematically
superfluous (anomaly \# 2). Basic questions about sets produced via the
power set axiom are apparently unanswerable (anomaly \# 1).

The best argument for the power set axiom, of course, is that it is needed
to go from $\Nb$ to $\Pc(\Nb)$, from which we can then construct the real line
$\Rb$, and if we accept it when applied to countable sets we presumably ought
to accept it when applied to any set. Now core mathematics obviously requires
the real line. However, it turns out that in the absence of the power set
axiom it is still possible to treat $\Pc(\Nb)$ and $\Rb$ as, in effect, proper
classes (see Sections \ref{sect4} and \ref{sect5}). One might imagine that
this approach would be inadequate for the development of normal mathematics,
but this is not correct. In fact, {\it when we eliminate the power set axiom
the result is a universe that matches actual mathematical practice with
remarkable accuracy} (see Section \ref{sect5}). As shocking as rejecting power
sets may seem at first, once one gets used to the idea it is highly satisfying
both philosophically and mathematically.

I address the philosophical justification (indeed, I would say, the
necessity) for doing this next.

\section{What are sets?}\label{sect2}

\begin{quote}
A pack of wolves, a bunch of grapes, or a flock of pigeons are all examples
of sets of things.

--- Paul Halmos (\cite{Hal}, p.\ 1)
\end{quote}

\begin{quote}
[I]t seems absurd to think of a collection as an entity distinct from the
items which compose it, so that when one buys a pair of shoes, one is buying
three things, the right shoe, the left shoe, and the pair. Not only that,
but if we do conceive of classes as entities and also make what would then
seem the reasonable assumption that they are to be included in the totality
of all things that there are, we fall into contradiction.

--- A.\ J.\ Ayer (\cite{Aye}, p.\ 43)
\end{quote}

\begin{quote}
Platonism is the medieval metaphysics of mathematics.

--- Solomon Feferman (\cite{Fef2}, p.\ 248)
\end{quote}

What are sets? The naive view expressed in the Halmos quote above is commonly
seen in elementary treatments of set theory. Very likely statements like this
are made for didactic reasons by authors who might otherwise adopt a more
nuanced position. On the other hand, this sort of assertion is so common that
one has to wonder how many mathematicians would simply accept it at face value.

It is false. A set of pigeons is {\it not} the same thing as a flock of
pigeons. Perhaps this is most obvious in the case of a set with exactly one
element. We would never speak of a flock consisting of a single pigeon and
if we did it would probably be synonymous with the pigeon itself. Yet a set
consisting of a single pigeon is, in mathematical usage, something totally
different from that pigeon. In fact it is not a physical entity at all, nor
is any set in the mathematical sense. We can say that a flock of pigeons is
overhead but we cannot say that a set of pigeons is overhead because this set
is, evidently, some sort of abstract entity which does not exist in time and
space. (Any reader who thinks sets are physical entities is invited to explain
the difference between a pigeon, the set whose sole member is that pigeon,
and the set whose sole member is the set whose sole member is that pigeon.)

In fact, as Ayer suggests, it seems a mistake to suppose that a collection,
as the word is used in ordinary speech, is any sort of {\it thing} at all,
distinct from the items which compose it. In ordinary speech it is generally
possible to reword any assertion concerning collections in an equivalent way
that avoids the concept entirely. For example, instead of saying that there is
a flock of pigeons overhead one could say: approximately thirty pigeons, in
close proximity to one another and flying in approximately the same direction,
are overhead. This suggests that the ``flock'' formulation is merely a
convenient abbreviation and assuming that the word refers literally to some
abstract metaphysical entity is an error on a par with, say, supposing that
``the average taxpayer'' is an actual person.

Thus, the founding myth of Cantorian set theory --- {\it that there are such
things as ``sets''} --- is false and absurd. This might make it look as though
the set concept has to be abandoned entirely, but there are ways of rescuing
it, at least partially. However, it does show that sets are not {\it canonical
objects}. When we refer to the set of all prime numbers, for instance, we are
not identifying a canonical abstract entity which exists in some metaphysical
arena. The mystical idea that the set of all prime numbers is some ideal
platonic {\it thing} is sometimes called ``realism'', but I do not see anything
realistic about it.  Defending it requires a double rationalization, first
to explain why it is not simply nonsensical, and second to explain why the
classical paradoxes do not discredit it. Perhaps when the power and economy
of the conceptualist view are better appreciated there will be less desire
to attempt such a defense.

I see two ways of handling the set concept rigorously. One is to replace
any expression that refers to a set by a (probably longer) expression
that makes no such reference, in the manner suggested earlier. This should be
easy to do in ordinary language settings but it is probably not workable in
the mathematical context where the set concept is used iteratively (i.e., we
consider sets of sets). In this case we must make the set concept rigorous
in another way, namely, by adopting a convention that when we say ``set'' we
really mean something else. For example, we might define the phrase ``a set
of natural numbers'' to mean ``an infinite sequence of 0's and 1's''. Under
this interpretation, the ``set'' of prime numbers would literally be the
sequence $001101010001\ldots$ and we could then carry out basic operations
involving ``sets'' of numbers without invoking any metaphysical platonic
universe of sets. (One could ask whether this particular proposal does not
merely trade one form of platonism for another. But this would only be the
case if one believed that the concept of a sequence of 0's and 1's cannot
be understood in other than platonic terms, which is surely false.)

The actual details of the convention that is adopted may seem to be of only
specialized interest, since the result is supposed to justify ordinary
mathematical usage and if it does so the working mathematician could safely
ignore it. However, this is not the case; there are important
consequences to taking the position that sets are not preexisting canonical
abstract entities but instead have to be defined or constructed in some way.
For example, we can no longer simply assert the ZFC axioms as being intuitively
true. We have to verify that it is possible to define or construct sets in
a way that actually satisfies these axioms. We also, {\it obviously}, have
to be sure that there is no circularity in our construction.

\section{Mathematical conceptualism}\label{sect3}

\begin{quote}
The basic idea is that all that matters is that the construction be
theoretically sound, so that one knows exactly what it would mean to
perform it and even to carry it to completion $\ldots$ All that matters
is clarity of ideas, and merely practical considerations are out of
place except insofar as they contribute to the ideas.

--- Lawrence Pozsgay (\cite{Poz}, p.\ 324)
\end{quote}

\begin{quote}
$\ldots$ the position that reason demands that we conform to the
vicious-circle principle, that nonconformity is unreasonable, illogical,
or scientifically unsound, and that the paradoxes are the result of failure
to adhere to this principle.

--- Charles Chihara (\cite{Chi}, p.\ xiv)
\end{quote}

Once we reject the idea that the world of sets is a well-defined independently
existing entity, and agree that any domain in which set-theoretic reasoning is
to take place must be in some sense constructed, we are then faced with the
question of just what sorts of putative constructions we should consider
legitimate. The basic premise of mathematical conceptualism is that in order
for a construction to be considered valid it need not be physically realizable,
but it must be {\it conceptually definite}, meaning that we must be able to
form a completely clear mental picture of how the construction would proceed.
This has nothing to do with any occult beliefs about mental states, but rather
is simply an expression of the idea that the proper domain of mathematics is
the realm of what is logically possible, and the proper criterion for
determining whether something is logically possible is not whether it can
be actually achieved in the physical universe (though this is certainly a
sufficient condition), but whether it can be conceived as taking place in
some concretely imagined world. What we are actually able to do in the
universe that we happen to inhabit is a merely contingent matter that we
should not expect to coincide with logical possiblity. On the other hand
we do insist on the absolute conceptual clarity of any imagined mathematical
construction.

The most straightforward way to do set theory on these terms would be to
explicitly describe some domain in a very concrete way, explicitly identify
for any two objects in the domain whether the first is to be considered an
element of the second, and then restrict all set-theoretic reasoning to take
place in that setting. (Let me say now that in this discussion I am going to
use the word ``domain'' synonymously with ``collection'' in the ordinary
language sense. As I discussed earlier, this is merely a convenient
abbreviation that could be eliminated throughout in favor of more complicated
but more literal expressions.) I describe one possible construction of this
type, $J_2$, in \cite{Wea1} and indicate there how core mathematics can be
developed entirely within it (see Section \ref{sect5}). Several other
more or less similar approaches are given in the references to \cite{Wea1}.
This shows that the rejection of a platonic set-theoretic universe does not
threaten ordinary mathematical practice.

Of course, in setting up a domain such as $J_2$ in which one intends to
carry out set-theoretic reasoning, one must ensure that both the domain and
the membership relation are well-defined, and in particular that there is
no vicious circle in their definition. One natural way to accomplish this is
to build up the domain in stages, in such a way that any ``set'' (i.e., any
object in the domain) only has as elements ``sets'' that already appeared at
earlier stages. This is an expression of the {\it vicious-circle principle}
championed by Poincar\'e and Russell, which is the central tenet of their
predicativist philosophy and which states in essence that every set should
be definable in terms of other sets that are logically prior to it.

A conceptualist would certainly accept the vicious-circle principle but he
would not regard it as fundamental. On the conceptualist view it
is merely a secondary consequence of the more basic fact that there is no
preexisting universe of sets and any domain in which set-theoretic reasoning
is to take place must be set up in a sensible manner. There are many ways in
which such a construction could fail and it seems pointless to try to identify
them all. However, circularity is certainly one important problem to avoid.
I do not understand the complicated explanations given by some workers in
mathematical foundations as to why there is really nothing wrong with
impredicative constructions.

I mentioned in the introduction that conceptualism is a modern version of
predicativism, and I can now explain that comment. Since conceptualism
incorporates the fundamental principle of predicativism, it seems fair to say
that the former is a version of the latter. However, I feel there are good
reasons for introducing the new term. First, as I just explained, there is a
difference in emphasis. Predicativism arose in response to the set-theoretic
paradoxes and this led to a preoccupation with banning circular definitions;
conceptualism accepts the predicativist conclusion but does not take it as
a central principle. Second,
there are several varieties of predicativism and conceptualism
corresponds to only one of these. The source of this variation is
the question of which principles of set construction are to be taken as
basic, and specifically whether infinitary constructions are accepted.
Historically, predicativists have generally accepted constructions of
length $\omega$, although Poincar\'e is a prominent
exception to this rule. But it is consistent with the vicious-circle
principle to reject any kind of infinite procedure, and this has led to
criticism based on the complaint that predicativism can take on different
meanings depending on which basic constructions are allowed. The vulgar
form of this objection alleges that one has no better reason to adopt
predicativism ``relative to the natural numbers'' than, say, predicativism
``relative to finite sets'' or ``relative to the real numbers'', and since
all cannot be true then all must be false. The conceptualist stance is not
subject to this criticism since we do have a compelling reason for preferring
the middle course, namely, that we have a conceptually clear idea of what
it would be like to carry out a construction of length $\omega$, whereas we
do not have a similarly definite idea of what it would be like to carry out
a construction of length $2^{\aleph_0}$. (Determining exactly which lengths
of constructions are legitimate on conceptualist grounds is rather subtle;
see \cite{Wea2}.)

Another important issue is the legitimacy of classical versus intuitionistic
logic. The predicativist literature tends to presume classical logic,
although there are occasional exceptions. This seems strange because we are
put in the position of affirming as definitely either true or false statements
that involve quantification over sets which we consider not well-defined.
(This kind of conceptual difficulty with quantification is especially acute
in \cite{Fef}, for example.) As I will indicate in the next section,
conceptualism entails use of a specific mixture of classical and intuitionistic
logic, with the former being applied to statements that refer only to sets
which are already within our purview and the latter to statements about sets
of some general type that might be constructed at some future time. This also
seems a good reason to adopt a special term for the conceptualist stance.

\section{Sets and classes}\label{sect4}

\begin{quote}
The nature of classes --- particularly proper classes, collections ``too
large'' to be sets --- is a perennial problem in the philosophy and
foundations of set theory $\ldots$ Existing theories of sets and classes
seem unsatisfactory because their `proper classes' are either indistinguishable
from extra layers of sets or mysterious entities in some perpetual, atemporal
process of becoming.

--- Penelope Maddy (\cite{Mad}, p.\ 299)
\end{quote}

Domains like $J_2$, while adequate for normal mathematics, may be
unnecessarily restrictive. On the one hand, they are satisfyingly definite
because in principle ``we have a complete and clear mental survey of all the
objects being considered, together with the basic interrelationships between
them'' (\cite{Fef}, p.\ 70). This comes out of the fact that we are actually
able to imagine {\it generating} $J_2$ one element at a time (in particular,
it is countable). On the other hand, if we insist on this condition then we
cannot include the power set of $\Nb$ as an element of the system because we
do not have ``a complete and clear mental survey'' of all sets of natural
numbers (say, realized as infinite sequences of 0's and 1's). But need we
insist on this? Even if we cannot concretely {\it generate} all possible
sequences, and therefore ``survey'' them, we can at least {\it unambiguously
recognize} any sequence presented to us for what it is. Is this enough to
justify conceptualistically regarding $\Pc(\Nb)$ as a set?

The preceding suggests that conceptualism may really support two distinct
set concepts. A set of the first type, which we may temporarily
call a set {\it in the strong sense}, would have the property that we can in
principle explicitly generate and survey its elements. In contrast, a set
{\it in the weak sense} would have the property that we can in principle
always unambiguously determine whether any given object belongs to the set.
This would certainly be the case if we could mentally survey its elements
as we could then simply check the given object against each element in turn;
however, as the power set of $\Nb$ demonstrates, that condition is not
necessary. I do not believe there is any reasonable weaker version of the
set concept that is consistent with our rejecting the idea of a preexisting
platonic universe of sets.

Please note that I am not asking whether it is possible to ``collect'' all
subsets of $\Nb$ into a ``totality'' (a question I consider meaningless,
though this sort of language is used in the predicativist literature).
Rather, the question is whether we may regard a set as having been specified
if we are able to recognize when something is an element but do not know how
to generate all of its elements, as in the case of the power set of $\Nb$.
This would not be an issue if we were working within a concretely specified
domain like $J_2$, since in a setting like this one could always in principle
generate all available subsets of $\Nb$ simply by working
through the entire domain one element at a time, determining of each element
in turn whether it is a subset of $\Nb$. The weak version of the set concept
would only be of interest in an open-ended setting which admitted indefinite
extension.

It seems to me that it is entirely up to us to decide whether we want to use
the word ``set'' in the weak or the strong sense. However, the weak sense is
problematic because some basic set-theoretic constructions fail for it. For
example, suppose we have a set $X$ in the weak sense whose elements are all
sets in the strong sense. Can we take the union of all $x \in X$? Evidently
not. A moment's thought shows that the result need not even be a set in
the weak sense.

A related example concerns least upper bounds in $\Rb$. Suppose we have a
subset of $\Rb$ which is a set only in the weak sense and which is known to
be bounded above. Regarding elements of $\Rb$ as lower Dedekind cuts, the
least upper bound of this subset is precisely its union, which we have just
seen is generally not a set even in the weak sense. Thus,
the least upper bound principle is valid only if we use the word ``set'' in
the strong sense, in which case $\Rb$ is not a set. This analysis, which goes
back to Poincar\'e, should begin to explain my comment in the introduction
about routine violations of the principle of building sets up from below.

Let me give a specific example of a subset of $\Nb$ which cannot be built
up from below because it apparently cannot be defined except in terms of
$\Pc(\Nb)$, which obviously is not a logically prior set in any sense.
Consider the language of second order arithmetic. This language has two kinds
of variables: some range over numbers ($\Nb$) and some range over sets of
numbers ($\Pc(\Nb)$). {\it Atomic formulas} are expressions of the form
$p_1 = p_2$, $p_1 \leq p_2$, or $p_1 \in X$, where $p_1$ and $p_2$ are
polynomials in the number variables with integer coefficients and $X$
is a set variable. All {\it formulas} (i.e., legal expressions) are then
built up from the atomic formulas using the logical symbols $\wedge$ (and),
$\vee$ (or), $\neg$ (not), $\Rightarrow$ (implies), $\Leftrightarrow$ (if
and only if), $\forall$ (for all), and $\exists$ (there exists). A
{\it sentence} is a formula in which all variables (of both types) are
quantified. We can enumerate the sentences by ordering them by length, and
lexicographically among those of equal length. Call this enumeration
$(\Ac_n)$. Finally, define
$$S = \{n \in \Nb: \Ac_n\hbox{ is true}\}.$$

One can begin to appreciate the complexity of $S$ if one realizes that all
normal mathematics can be encoded in second order arithmetic \cite{Sim2}. Thus,
$S$ contains the answer to Goldbach's conjecture, the Riemann hypothesis,
and whether $P = NP$; the solution to the word problem for all finitely
presented groups; an index of which Turing machines halt on which inputs;
and so on. As a matter of fact, these are all first order questions, i.e.,
they can be formulated without using set variables, and the second order
language is incomparably stronger. It is so strong, in fact, that the
definition of $S$ is fundamentally circular. $S$ is ``impredicative'' because
{\it in order to diagnose whether a given formula $\Ac_n$ is true we need
to quantify over all sets of numbers, which requires that we already have
access to $\Pc(\Nb)$}. To put this another way, if we assume that there is such
a thing as the set of all subsets of $\Nb$, this would imply that $S$ does
exist and hence that there are subsets of $\Nb$ which cannot be defined
except by reference to $\Pc(\Nb)$. So this assumption involves a genuine
circularity. From the conceptualist standpoint, we cannot form any precise
mental picture of how $S$ could be constructed, because any construction
would at some point require that we already possess all subsets of $\Nb$.
Thus, {\it we cannot build $\Pc(\Nb)$ up from below}.

It has been suggested that $\Nb$ is just as impredicative as $\Pc(\Nb)$ since
it is possible to define a natural number in terms of the set of all natural
numbers. In other words, one can give a definition of some $n \in \Nb$ such
that in principle we need to search through all of $\Nb$ in order to
determine the value of $n$, just as we needed to use $\Pc(\Nb)$ in order
to identify the set $S$. But this misses the point because whatever the
value of $n$ is, it {\it could} be built up from below in a trivial way.
Once we agree that this is so for every natural number, we can build up
$\Nb$ from below. And once $\Nb$ is available, we can legitimately use
impredicative definitions to identify particular natural numbers. The
example involving the true sentences of second order arithmetic is not
analogous because there is {\it no alternative definition} of $S$ that
avoids circularity. It does not help to say that $S$ could
be built up from below by simply enumerating its elements, since this
assumes that we already know all possible ways of enumerating subsets
of $\Nb$, which amounts to assuming that $\Pc(\Nb)$ is already available.

Now, returning to the dilemma posed earlier about the two possible
conceptualistic set concepts, an appealing solution is to reserve the word
{\it set} only for sets in the strong sense, and to call sets in the weak
sense {\it classes} (and classes which are not sets {\it proper classes}).
Then the power set of $\Nb$ will be a proper class, as will the real line,
and indeed, apparently, any object which is classically uncountable, provided
that object exists conceptualistically at all. I find this crisp distinction
between the conceptualist notions of ``set'' and ``class'' quite satisfying
in comparison with their hazy classical interpretation.

A further question that could be asked is whether it is possible in a
conceptualistic system to use variables which range over proper classes such
as $\Pc(\Nb)$. I consider this question in \cite{Wea2} and reach the
conclusion that this is indeed possible, but only if intuitionistic logic
is used. If we are working in an open-ended domain that is indefinitely
extendible then the law of the excluded middle is generally not justified.

In sum, there are two essentially different ways of straightforwardly setting
up conceptualistic set-theoretic systems. The first approach involves reasoning
about a concretely described closed domain of which one can in principle make
a complete mental survey (``in principle'' since one need not actually make
such a survey, one need only be able to concretely imagine making it);
in such a setting one can reason using classical logic. The example
I have been using in this case is the domain $J_2$ as described in \cite{Wea1}.
The second approach involves reasoning about an open-ended domain with the
property that one can recognize putative elements but one cannot generate or
survey them all. In this case intuitionistic logic must be used (although
classical logic is legitimate to some extent; see \S 2.2 of \cite{Wea2}).
A basic example of this type is the system ${\rm ACA}_0$ (see \cite{Sim2}),
using intuitionistic logic supplemented by (1) $\Ac \vee \neg \Ac$ for all
atomic formulas $\Ac$ and (2) the ``numerical omniscience schema''
$$(\forall n)(\Ac(n) \vee \neg \Ac(n)) \Rightarrow
[(\forall n)\Ac(n) \vee (\exists n)\neg \Ac(n)]$$
for all formulas $\Ac$ of second order arithmetic and all number variables
$n$. Some stronger conceptualistic systems using intuitionistic logic are
discussed in \cite{Wea2}.

\section{Core mathematics}\label{sect5}

\begin{quote}
The main purpose of a construction and development of the theory will
be not so much the exhibition of a formal system as the basis of all
mathematics, as the presentation of an argument to justify all mathematical
reasoning which does not get into the transcendent $\ldots$ For example,
as the Bourbaki group continues to turn out more and more volumes of their
treatise, we show for each volume how all the definitions and proofs can
be formalized in the theory $\Sigma$.

--- Hao Wang (\cite{Wan}, p.\ 582)
\end{quote}

At this point the reader may be curious as to how it is possible to develop
core mathematics within a conceptualistic system. In the last section
I indicated two ways of going about this. (A third possibility will be
discussed in Section \ref{sect6}.)
One is to work in a fixed domain using classical logic and the other
is to work in an open-ended domain using intuitionistic logic. There
is a good deal to be said for the second alternative; one can make a
case that this is the more elegant approach. However, the use of
intuitionistic logic, besides being relatively unfamiliar to most
mathematicians, does complicate matters somewhat. So I will outline
here a development of ordinary mathematics using the first method. For
more details the reader is referred to \cite{Wea1}.

We start with the set $J_1$ of all hereditarily finite sets. That is,
$$J_1 = \emptyset \cup \Pc(\emptyset) \cup \Pc(\Pc(\emptyset)) \cup \cdots.$$
Sets in $J_1$ can be associated to finite trees in the following way:
given a finite tree, label each terminal node with the empty set, and
inductively label every node with the set of labels of its immediate
successors. Then the root node (indeed, every node) will be labelled
with a hereditarily finite set, and every hereditarily finite set can
be obtained from a finite tree in this way.

The conceptualistic mini-universe $J_2$ in which we propose to develop core
mathematics is obtained from $J_1$ by applying certain closure operations.
We begin by placing every hereditarily finite set in $J_2$, as well as the
set $J_1$ itself. The set $J_2$ is then generated from this foundation by
repeated application of a small number of ``rudimentary'' operations; for
instance, if $x$ and $y$ are in $J_2$ then so are $\{x,y\}$, $x - y$,
$x \times y$, and so on. Two remarkable facts emerge: first, a surprisingly
large variety of set-theoretic constructions can be built up from the
primitive operations specified, and second, the resulting set $J_2$ is
closed under {\it first-order definability}. Loosely speaking, this means
that if a set $x$ belongs to $J_2$ then every subset of $x$ that we can
explicitly name also belongs to $J_2$. As a result, even though $J_2$ is
countable, most ordinary mathematical constructions will not lead one outside
of it.

Since $J_2$ is so small, I use the prefix $\iota$ (iota) to identify the
$J_2$-versions of classical mathematical notions. For example, an
{\it $\iota$-set} is just an element of $J_2$. An {\it $\iota$-class}
is a first-order definable subset of $J_2$ whose intersection with every
$\iota$-set is an $\iota$-set. (The last condition is needed because,
roughly speaking, a definable subset of $J_2$ can fail to be an element
of $J_2$ either because of complexity or because of size. We want the
$\iota$-classes to be just those which are too large to be $\iota$-sets.)
The canonical example of a proper $\iota$-class, i.e., an $\iota$-class
which is not an $\iota$-set, is $\Pc_\iota(\Nb) = \Pc(\Nb) \cap J_2$. Any
$\iota$-set of $\iota$-subsets of $\Nb$ could be diagonalized using
rudimentary functions to produce a new $\iota$-subset of $\Nb$, which
shows that $\Pc_\iota(\Nb)$ cannot be an $\iota$-set.

The preceding argument makes use of the fact that every $\iota$-set is
$\iota$-countable, i.e., there exists a bijection in $J_2$ between $\Nb$
and any given $\iota$-set. This illustrates the regularity of $J_2$.
Another such illustration is the existence of an $\iota$-class which
well-orders $J_2$. In other words, a strong form of the axiom of choice
holds in $J_2$ according to which not only each set, but the universal
class, can be well-ordered.

The extent to which we can do mathematics within $J_2$ depends crucially on
our ability to manipulate $\iota$-classes. Just as with proper classes in the
classical sense, there are limits to what one can legitimately do with
$\iota$-classes. However, basic constructions are possible: if $X$ and $Y$
are $\iota$-classes then so are $X \cup Y$, $X \cap Y$, $J_2 - X$,
$X \times Y$, and the ``power class''
$$\Pc_\iota(X) = \{x \in J_2: x \subseteq X\}.$$
Moreover, if $x$ is an $\iota$-set and $X \subseteq x \times J_2$ is an
$\iota$-class then $\bigcup X_a$ and $\bigcap X_a$ are $\iota$-classes,
where $X_a = \{z: \langle a,z \rangle \in X\}$, as is
\begin{eqnarray*}
&&{\prod}_\iota X_a = \{f \in J_2:
f\hbox{ is a function with domain $x$}\\
&&\qquad\qquad\qquad\qquad\qquad\qquad
\hbox{such that $f(a) \in X_a$ for all }
a \in x\}.
\end{eqnarray*}

It is crucial in the last comment that $x$ be an $\iota$-set and not
a proper $\iota$-class. Informally, we could say that the union,
intersection, and $\iota$-product of an $\iota$-set of $\iota$-classes is
again an $\iota$-class. This statement is not literally accurate
because the $X_a$ could be proper $\iota$-classes and there is no
such thing as an $\iota$-set of proper $\iota$-classes. However, each
$\iota$-class $X_a$ is indexed by an element $a \in x$; we call
$a$ a {\it proxy} for $X_a$ and interpret a statment like ``the
union of an $\iota$-set of $\iota$-classes is an $\iota$-class''
as meaning ``given an $\iota$-set of proxies for $\iota$-classes,
the union of the corresponding $\iota$-classes is an $\iota$-class''.

Direct constructions establish that $\Nb$, $\Zb$, and $\Qb$ are
$\iota$-sets. We define $\Rb_\iota$ to be the family of Dedekind cuts in
$J_2$; this is a proper $\iota$-class.

What other basic mathematical objects exist as $\iota$-sets and
$\iota$-classes? Any finite group, and most if not all countable groups of
interest, are $\iota$-sets. Indeed, most if not all countable structures of
any sort that appear in mainstream mathematics are $\iota$-sets. Since we can
take finite products of $\iota$-classes, $\Rb_\iota^n$ is an $\iota$-class
for any $n \in \Nb$, and because of the good closure properties of $J_2$
under definability, any variety in $\Rb_\iota^n$ defined by a real polynomial
with coefficients in $\Rb_\iota$ is an $\iota$-class. Because we can take
countable $\iota$-products, $\Rb_\iota^\Nb$ is also an $\iota$-class,
and it contains $\iota$-subclasses $l^p_\iota$ for $1 \leq p \leq \infty$.
Moreover, even though we cannot take an $\iota$-product of a proper
$\iota$-class of $\iota$-classes, we can still define $l^p_\iota(X)$ for
any $\iota$-class $X$ and any $1 \leq p < \infty$. Namely, we let it be the
set of functions in $J_2$ whose domain is an $\iota$-subset of $X$, whose
range is contained in $\Rb_\iota$, and which is $p$-summable; this turns out
to be an $\iota$-class. This construction does not work for $l^\infty(X)$,
however, since elements of this space cannot be assumed to have countable
support.

What about morphisms? We define an {\it $\iota$-function} between two
$\iota$-classes to be a function whose graph is an $\iota$-class and which
satisfies the regularity property that the image of any $\iota$-set is an
$\iota$-set. We can then define structures such as $\iota$-metrics. If $X$
is an $\iota$-class, then an {\it $\iota$-metric} on $X$ is (essentially;
see \cite{Wea1} for details) an $\iota$-function
$D: X \times X \to \Rb_\iota^+$ which satisfies the usual
metric axioms. We can carry out a Cauchy sequence construction to complete
any $\iota$-metric space, although actually this construction involves
Cauchy {\it $\iota$-sequences}, i.e., Cauchy sequences which are
$\iota$-functions from $\Nb$ into $X$.

The theory of $\iota$-metric spaces is generally similar to the theory of
classical metric spaces. However, proofs are sometimes a little trickier and
hypotheses sometimes must be strengthened, owing to the rigors of reasoning
within $J_2$. For example, classically we know that any closed subset of a
separable metric space is separable; in $J_2$ the corresponding assertion
about closed $\iota$-subclasses is not so obvious because we have to
construct the desired dense {\it $\iota$-subset} in terms of a closed
{\it $\iota$-subclass}. However, it turns out that this can be done if
we assume that the ambient metric space is boundedly compact, i.e., every
closed ball is compact. Now it is
certainly possible that a mainstream mathematician might want to use this
theorem without the assumption of bounded compactness. But it seems that
results that do not hold in $J_2$ are in fact actually used only relatively
rarely. Moreover, it seems that {\it really necessary} uses are vanishingly
rare in mainstream mathematics. This can only be definitively established by
actually translating large sections of core mathematics into $J_2$, which has
not yet been done. (Of course, to some extent it also depends on one's
definition of ``mainstream'' mathematics.)

More abstract structures can also be defined and studied within $J_2$. For
example, using the proxy technique we can define $\iota$-topologies. An
{\it $\iota$-topology} on an $\iota$-class $X$ will be an $\iota$-subclass
$\Tc$ of $T \times X$ for some $\iota$-class $T$. We think of the elements
of $T$ as proxies for the {\it $\iota$-open} $\iota$-subclasses of $X$, which
are defined to be the sections of $\Tc$. That is, for each $a \in T$ we have an
$\iota$-open $\iota$-class $U_a = \{t \in X: \langle a,t \rangle \in \Tc\}$.
We require that $\emptyset$ and $X$ be $\iota$-open, that any union of an
$\iota$-set of $\iota$-open $\iota$-classes be $\iota$-open, and that any
intersection of finitely many $\iota$-open $\iota$-classes be $\iota$-open.
A generally similar approach allows us to define $\iota$-$\sigma$-algebras
and $\iota$-measure spaces.

The reader should now be getting some sense of how core mathematics can
be developed within a conceptualistic setting like $J_2$. I want to mention
one other phenomenon that I find rather compelling. In any conceptualistic
system in which the real line is treated as a proper class, the same will
clearly be true of any (nonzero) Banach space $E$. This presents an
apparent problem
for duality theory, since the elements of the dual of $E$ are continuous
linear functions from $E$ into the scalar field (say $\Rb_\iota$), and if $E$
is a proper class then the graph of any such function will be a proper class
and we will not be able to collect them all into a single class $E'$. Thus,
because proper classes cannot be elements of other classes, in general we
cannot form the dual space, although we can handle individual elements of
the dual space. However, this problem can be effectively resolved if $E$ is
separable as in that case any continuous function will be determined by
its values on a countable subset $E_0$, so we can take as proxies for the
continuous linear functionals their restrictions to $E_0$; these restrictions
will be {\it sets} and we can legitimately form the class of all such
restrictions. Thus, we have a workable construction of the dual space of
any {\it separable} Banach space.

Since the dual of a separable Banach space may or may not be separable,
in general we cannot form its second dual. This also seems like a problem
because double dual techniques are quite standard in functional analysis.
However, close examination of typical double dual arguments reveals that
{\it they generally involve only individual elements of the second dual,
not the entire second dual as a Banach space}. Thus, our limited ability
to conceptualistically handle duals of non-separable Banach spaces generally
{\it is} sufficient to support typical double dual arguments provided the
original Banach space is separable. For example, the most important result
about second duals states that the unit ball of $E$ is weak* dense in the unit
ball of $E''$. This could be rephrased as: any unit norm continuous linear
functional on $E'$ can be weak* approximated by unit norm elements of $E$ ---
and this is conceptualistically expressible and true (for example, in $J_2$).

One could argue that the classical picture is simpler and more elegant:
classically, every Banach space has a dual, so that starting from $E$ one
can form $E'$, $E''$, $E'''$, etc. But what I find impressive is that in
practice arguments which use duals higher than the second are quite rare, and
that the conceptualistic universe cuts off at {\it precisely} the point ---
partway into the second dual --- where typical mainstream usage cuts off.
This is an example of what I meant when I said in the introduction that the
conceptualistic universe matches mathematical practice with remarkable
accuracy.

I must emphasize that $J_2$ is not the only possible conceptualistically valid
system. It is just an example. There are many possibilities and a decision
as to which really is the best suited to ordinary mathematics might be
premature at this point. The main advantage of $J_2$ is supposed to be its
retention of the standard language of set theory. Type-theoretic systems,
such as $W$ of \cite{Fef2}, seem further from the classical intuition, while
work in subsystems of second order arithmetic \cite{Sim2}  involves a more
elaborate coding machinery. (Our use of proxies is also a form of coding, but
it is less pervasive: compare the concept ``continuous function between
complete metric spaces'' in the two approaches. One could also argue that
defining a measure to be a linear functional on $C(X)$, as one does in
other approaches, represents a loss of content. This should be a theorem,
not a definition.)

\section{Conceptualistic reductionism}\label{sect6}

\begin{quote}
It thus became apparent that the ``finite Standpunkt'' is not the only
alternative to classical ways of reasoning and it is not necessarily
implied by the idea of proof theory. An enlarging of the methods of
proof theory was therefore suggested: instead of a restriction to
finitist methods of reasoning, it was required only that the arguments
be of a constructive character $\ldots$.

--- Paul Bernays (\cite{Ber}, p. 502)
\end{quote}

In the introduction I suggested that there is an analogy between the
transition from Cantorian set theory to conceptualistic mathematics and the
transition from a naive use of infinitesimals to the epsilon-delta method.
This analogy is surprisingly robust. The concept of infinitesimals is grounded
in a real process, allowing a variable to decrease to zero, which we imagine
in some vague way as extending to an ``idealized'' point at which the variable
becomes infinitesimal (but not zero). Set-theoretic objects like $\Pc(\Nb)$
and $\aleph_1$ can be understood in similar terms as ``idealized'' limits of
legitimate processes which do not in fact have limits. Given any positive
real number, we can always divide by two to get a smaller positive real
number; given any conceptualistically legitimate set of countable ordinals,
we can add one to its union to get a larger countable ordinal. Classically,
in some vague way we imagine the latter process as extending to an idealized
point at which the set of countable ordinals becomes uncountable (but not a
proper class), but conceptualistically we do not recognize this limit as
meaningful.

The diagnosis is analogous and the treatment is also analogous. A rigorous
interpretation of reasoning using infinitesimals replaces them with variable
quantities which can always be made smaller, and adds precision by formulating
the relationships between small positive values in epsilon-delta terms.
Conceptualistically we work in the opposite direction and formulate the
possible techniques of set construction in precise terms, giving rise to a
conception of a set-theoretic universe whose growth is stratified in a far
finer manner than the cumulative hierarchy of ZFC. The epsilon-delta analog
can be seen in conceptualistic verifications of statements of the form
$(\forall X)(\exists Y)\Ac(X,Y)$, which classically take place in some
idealized fully formed universe but conceptualistically are made rigorous by
specifying, for a given set $X$ which appears at some level of the hierarchy
(``for every $\epsilon$''), how much farther up in the hierarchy we must go
(``there exists a $\delta$'') in order to find a set $Y$ which satisfies
$\Ac(X,Y)$.

But there is another way of making infinitesimals rigorous: non-standard
analysis. In this approach we rigorously construct a legitimate model of a
``non-standard real line'' which actually contains infinite and infinitesimal
quantities in addition to ordinary real numbers. The intricacy of the naive
notions is reflected in the sophistication of the construction, which uses
powerful tools to produce a rather elaborate domain in which the naive
arguments work. The conceptualistic analog of this approach is the following.
Starting with some classical axiom system such as ZFC, we construct in a
conceptualistically legitimate manner a ``non-standard'' model of the axioms.
This, in a way, justifies naive reasoning in the original system, which is now
understood as conceptualistically valid reasoning within an elaborate
non-standard domain.

Now at present there seems little prospect of conceptualistically validating
reasoning within ZFC in such a manner. However, through the work of
Friedman, Simpson, and others in the field of ``reverse mathematics''
\cite{Sim2}, we know that classical reasoning generally does not require
the full strength of ZFC but can in fact be carried out in substantially
weaker systems. In fact, the goal of reverse mathematics is to exactly
determine which theorems require which axioms. Of the five axiom systems
identified by Friedman as fundamental, three are conceptualistically legitimate
and two (${\rm ATR_0}$ and $\Pi^1_1-{\rm CA}_0$) are not. However, here there
is a real prospect of conceptualistically proving consistency by constructing
non-standard models. In the case of ${\rm ATR}_0$ this can be accomplished by
combining the construction of a non-standard model \cite{FAS} with an ordinal
analysis that conceptualistically legitimates the construction \cite{Wea2}. I
believe a conceptualistic consistency proof of $\Pi^1_1-{\rm CA}_0$ is also
within reach. This general program could be termed ``conceptualistic
reductionism''; it is a version of Hilbert's idea (to which Bernays refers
in the above quote) of finding a finitistic consistency proof of classical
analysis, in which ``finitistic'' is replaced by ``conceptualistic''.

The motivation for carrying out this conceptualistic version of Hilbert's
program is that it would, in a sense, justify classical theorems
which are not conceptualistically valid. As I indicated earlier, the vast
bulk of core mathematics can be fairly straightforwardly developed within
conceptualistically legitimate systems. However, there are some classical
theorems which cannot, probably the best known being the
Cantor-Bendixson theorem which asserts that every closed subset of $\Rb$
is the union of a countable set and a perfect set. (I hasten to add that
some previous work on classical theorems which are allegedly not provable
in predicatively legitimate systems is incorrect. See \cite{Wea2}.) My
feeling is that results of this type are dispensible and that core
mathematics can be developed in a fully satisfactory manner without them.
But another approach is to work in an axiom system that supports the
desired classical theorems and which can be conceptualistically shown to be
consistent. Perhaps there will come a time when many mathematicians will
regard such systems in the same way that they now regard the nonstandard reals.

Of course, no one is barred from working in systems like ZFC which are not
likely to ever be conceptualistically proven consistent. In these cases we
would just have to accept that we cannot be certain the system is consistent
and that we lack the really robust understanding of the system that would come
with possession of a concrete model.

This conceptualistic reinterpretation of Hilbert's program contributes to a
feature of conceptualism that I find quite appealing: its synthesis of
the three great twentieth-century foundational philosophies of platonism,
intuitionism, and formalism. In fact conceptualism embodies aspects of
all three. First, its approach to elementary number theory is essentially
{\it platonic}. Although we do not believe in the set of natural numbers as an
ideal platonic entity which exists in some metaphysical realm, we do regard
the structure of $\Nb$ as completely definite, so that all elementary
assertions about natural numbers are unambiguously either true or false.
On the other hand, when we pass to arbitrary sets we regard the universe
of discourse as inherently indefinite and always capable of extension, which
makes us {\it intuitionists} with regard to set theory. Finally, we
can justify classical axiom systems by proving them consistent, which makes
us {\it formalists}. I find it satisfying and appropriate that conceptualism
should in this way be a reinterpretation and development, and not a mere
rejection, of the philosophies of the past.

%##

\end{document}